\documentclass[11pt]{amsart}
\usepackage{amssymb,amsmath,amsopn,textcomp, thmtools}
\usepackage[foot]{amsaddr}
\usepackage[version=4]{mhchem}
\usepackage{verbatim}
\newtheorem{definition}{Definition}[section]
\newtheorem{theorem}{Theorem}[section]

\newcommand{\iR}{\mathbb{R}}
\newcommand{\iN}{\mathbb{N}}
\newcommand{\iZ}{\mathbb{Z}}

\newcommand{\cP}{\mathcal{P}}

\newcommand{\cI}{\mathcal{I}}

\usepackage{bbm}
\newcommand{\R}{\mathbb{R}}

\newcommand{\cH}{\mathcal{H}}

\usepackage{extarrows}
\usepackage{dsfont}
\usepackage{tikz,subcaption} 

\usepackage{mathabx} 

\allowdisplaybreaks

\begin{document}

\title{Li--Yau and Harnack estimates for nonlocal diffusion problems}
\author{Rico Zacher}
\thanks{The author gratefully acknowledges the hospitality of the Ghent Analysis \& PDE Center during his visit. }
\email[Corresponding author:]{rico.zacher@uni-ulm.de}
\address{Institut f\"ur Angewandte Analysis, Universit\"at Ulm, Helmholtzstra\ss{}e 18, 89081 Ulm, Germany.}


\begin{abstract}
These notes give a brief introduction to differential Harnack inequalities and summarise the main results of the mini-course “Li–Yau and Harnack estimates for nonlocal diffusion problems”, presented by the author at the Seasonal School on PDEs “Oscillation Phenomena, PDEs, and Applications: A Comprehensive School in Mathematical Analysis”, held at Ghent University in October 2025.
\end{abstract}

\maketitle

 

\section{Introduction: Classical Li--Yau and CD inequality}
We begin by recalling the classical Li--Yau inequality for the heat equation in $\iR^d$. 
Suppose $u:\,[0,\infty)\times \iR^d \to (0,\infty)$ is sufficiently smooth and solves $\partial_t u-\Delta u=0$ on $(0,\infty)\times \iR^d$. Then
\begin{equation} \label{CL1}
- \Delta (\log u) \leq \frac{d}{2t}\quad  \mbox{in}\;(0,\infty)\times \iR^d.
\end{equation}
Since $\partial_t(\log u)-\Delta(\log u)=|\nabla(\log u)|^2$, this is equivalent to
\begin{equation} \label{CLY}
\partial_t(\log u)\ge |\nabla (\log u)|^2-\frac{d}{2t}\quad \mbox{in}\;(0,\infty)\times \iR^d.
\end{equation}
These inequalities extend to complete $d$-dimensional Riemannian manifolds $M$ with nonnegative Ricci curvature, 
see the celebrated work by Li and Yau \cite{LY}. Note that (\ref{CL1}) is sharp, one has equality for the heat kernel $u(t,x)=(4\pi t)^{-d/2} \exp\big(\frac{-|x|^2}{4t}\big)$. \eqref{CLY} is referred to as parabolic gradient estimate or differential Harnack inequality. In fact, integrating 
\eqref{CLY} over the line connecting $(t_1,x_1)$ and $(t_2,x_2)$ with $0<t_1<t_2$ yields
the sharp Harnack estimate
\[
u(t_1,x_1)\le u(t_2,x_2) \Big(\frac{t_2}{t_1}\Big)^{d/2} \exp\Big(\frac{|x_1-x_2|^2}{4(t_2-t_1)}\Big).
\]

Let us give a sketch of a proof of \eqref{CL1}. As already mentioned, $v:=\log u$ solves
\begin{equation} \label{CLY1}
\partial_t v-\Delta v=|\nabla v|^2,
\end{equation}
by the chain rule $\Delta H(u)=H'(u)\Delta u+H''(u)|\nabla u|^2$ with $H=\log $. To see that
 $-\Delta v\le \frac{d}{2t}$, we apply $\Delta$ to \eqref{CLY1} and use Bochner's identity.
\begin{align*}
\partial_t \Delta v  -\Delta(\Delta v) & =\Delta\big(|\nabla v|^2\big)\\
& = 2\nabla v\cdot \nabla\Delta v+2|\nabla^2 v|_{HS}^2\;\quad\Big(+2\text{Ric}(\nabla v,\nabla v)\Big),
\end{align*}
where $|\cdot|_{HS}$ denotes the Hilbert-Schmidt norm and the term in brackets indicates the
situation of a manifold with an additional curvature term, which can be dropped if the manifold is nonnegatively curved. Now, $|\nabla^2 v|_{HS}^2\ge \frac{1}{d}\,(\Delta v)^2$ (this is an example of a CD inequality, see below), and thus
\[
\partial_t \Delta v -\Delta(\Delta v) \ge  2\nabla v\cdot \nabla\Delta v+\frac{2}{d}\,(\Delta v)^2,
\]
which shows that $\Delta v$ is a supersolution of a semilinear parabolic equation. The function
$\omega(t):=-\frac{d}{2t}$ solves $\partial_t \omega=\frac{2}{d}\omega^2$ and is therefore
a solution of this semilinear PDE. Applying the comparison principle then yields
$\Delta v\ge \omega$. Although the final step is only formal, it can be made rigorous by employing suitable cut-off functions.

We next explain the $\Gamma$-calculus and the notion of curvature-dimension (CD) inequality due to Bakry and \'Emery (\cite{BGL,BE}). Let $L$ be the generator of a Markov semigroup. Then, for $f$ and $g$ from a suitable class of functions, the carr\'e du champ operator and iterated carr\'e du champ operator associated with $L$ are defined by
\begin{align*}
\Gamma(f,g) &=\frac{1}{2}\big(L(fg)-fLg-gLf\big),\\
\Gamma_2(f,g) & =\frac{1}{2}\big(L\Gamma(f,g)-\Gamma(f,Lg)-\Gamma(g,Lf)\big).
\end{align*}
We also set $\Gamma(f)=\Gamma(f,f)$ and $\Gamma_2(f)=\Gamma_2(f,f)$. $L$ is said to satisfy the curvature-dimension inequality $CD(\kappa,d)$ with curvature parameter $\kappa\in\iR$ and dimension parameter $d\in [1,\infty]$ if
(here $\mu$ is a fixed invariant and reversible measure) 
\begin{equation} \label{CDI}
\Gamma_2(f)\ge \kappa \Gamma(f)+\frac{1}{d}(Lf)^2,\quad \mu-a.e.
\end{equation}
For example, $L=\Delta$ on $\iR^d$ satisfies $CD(0,d)$. Here, $\Gamma(f,g)=\nabla f \cdot \nabla g$ and $\Gamma_2(f)=|\nabla^2 f|_{HS}^2$.

Assuming in addition that $L$ is a diffusion generator in the sense that the chain rule
$L H(f)=H'(f)L f+H''(f)\Gamma(f)$ holds for all $H\in C^2$ and suitable $f$, one can show the following: (i) $CD(0,d)$ with $d<\infty$ implies the Li--Yau inequality $-L (\log e^{Lt} f)\le \frac{d}{2t}$ for positive $f$, and (ii) $CD(\kappa,\infty)$ with $\kappa>0$ (i.e.\ we have positive curvature) implies important functional inequalities such as the log-Sobolev and Poincar\'e inequality,
which in turn entail exponential decay of entropy and variance as well as hypercontractivity. We refer to
the excellent monograph \cite{BGL}.
\section{Nonlocal operators: examples and difficulties}
We now turn to nonlocal operators. Let $X$ be a metric space and consider a Markov generator
of the form
\[
Lf(x)=\int_X \big(f(y)-f(x)\big)\,k(x,dy),\quad x\in X.
\]
The essential feature here is the difference structure under the integral. Many important examples of nonlocal operators are of this form, for example
\begin{itemize}
\item [(i)] the fractional Laplacian: $X=\iR^d$, $\beta\in (0,2)$,
\begin{equation} \label{frakLaplacianDef}
Lf(x)= - \big(-\Delta)^\frac{\beta}{2}f(x) = c_{\beta,d}\,\, \mathrm{p.v.}\int_{\R^d} \frac{f(y)-f(x)}{|x-y|^{d+\beta}}\,dy,
\end{equation}
\item [(ii)] generators of Markov chains (generalised Laplacians on a graph):
$X$ is a finite or countably infinite set, 
\begin{equation} \label{chain}
Lf(x)=\sum_{y\in X}k(x,y)\big(f(y)-f(x)\big).
\end{equation}
\end{itemize}

Suppose that $u$ is a positive solution of 
\begin{equation} \label{heatL}
\partial_t u-Lu=0, \quad t>0.
\end{equation}
Under which conditions does a Li--Yau inequality hold? Moreover, if such an estimate is available, does it imply a parabolic Harnack inequality?
For the fractional Laplacian $L= - \big(-\Delta)^\frac{\beta}{2}$, the question of the validity of a 
Li--Yau inequality was highlighted as a major open problem in a survey by Garofalo \cite{Garo}.

In the nonlocal case, several substantial difficulties arise. A major obstacle is the failure of the chain rule for the operator $L$. A second issue concerns the correct formulation of a Li--Yau type inequality in this context. Moreover, the classical curvature-dimension condition  $CD(0,d)$
turns out to be insufficient; this phenomenon was already observed in the graphical setting in \cite{Harvard}, indicating the need for stronger or modified curvature-dimension conditions. Finally, in order to derive a Harnack inequality, one must work with a kind of ``nonlocal gradient''.

The first positive results in this direction were obtained for finite and locally finite graphs about a decade ago. Yau and several co-authors (\cite{Harvard}) circumvent the failure of the chain rule by working with the square root of $u$ instead of $\log u$. They introduce the so-called exponential curvature-dimension conditions $CDE(\kappa,d)$ and $CDE'(\kappa,d)$. Under these conditions,
they establish Li--Yau inequalities with a relaxation term of the form $\frac{C}{t}$ for a wide class of graphs, and are moreover able to derive a parabolic Harnack inequality. Later, M\"unch introduces 
a general $\Gamma^\psi$-calculus for concave functions $\psi$ and an associated CD-condition $CD_\psi(\kappa,d)$. For $\psi=\sqrt{}$, he recovers some of the results from \cite{Harvard}, while the choice
$\psi=\log$ leads to new logarithmic Li--Yau inequalities. He also derives Harnack inequalities from the Li--Yau estimates. All of his results are restricted to finite graphs. Subsequently, Dier, Kassmann and the author of these notes (\cite{DKZ}) go one step further by establishing logarithmic Li--Yau inequalities based on a new CD-condition 
$CD(F;0)$, which involves a so-called CD-function $F$. This framework allows for more general relaxation functions (see below) and yields improvements over previous results; some of the resulting estimates are even sharp. As we will see, in the case of long-range jump operators, the flexibility provided by CD-functions proves to be crucial.

It should be noted that, concerning the fundamental problem of defining lower curvature bounds in the discrete setting, alternative approaches based on optimal transport have been developed (see, e.g., \cite{EM, Mie, Olli}).

We now turn to the notion of a CD-function. The basic idea is to replace the quadratic
dimension term in the CD inequality (see \eqref{CDI}) by a more general term. A continuous function $F:\,[0,\infty)\rightarrow [0,\infty)$ is called CD-function if $F(0)=0$, $F(x)/x$ is strictly
increasing on $(0,\infty)$,
and $1/F$ is integrable at $\infty$. A simple example is $F(x)=\nu x^2$ with $\nu>0$.

Given a CD-function $F$, there exists a 
unique strictly positive solution $\varphi$ of the ODE 
\begin{align*}
\dot{\varphi}(t)+F(\varphi(t))=0,\quad t>0,
\end{align*}
which has $(0,\infty)$ as its maximal interval of existence. This function 
$\varphi$ is strictly decreasing and log-convex, and it satisfies 
$\varphi(0+)=\infty$ and $\varphi(\infty)=0$, see \cite{DKZ}.  
The function $\varphi$ is called the relaxation function associated with $F$.

Following \cite{DKZ}, we next explain how the equation for the logarithm of a positive
solution $u$ of \eqref{heatL}
can be properly formulated in the nonlocal setting. Let us consider the Markov chain setting, that is,
$L$ is given by \eqref{chain}. An elementary calculation shows that $v:=\log u$ solves
\begin{equation} \label{loguEqu}
\partial_t v-L v= \Psi_\Upsilon(v)\quad \mbox{in}\;(0,\infty)\times X.
\end{equation}
Here, for $H:\,\iR\to \iR$,
\begin{equation} \label{PsiH}
\Psi_H (f)(x)=\sum_{y\in X}k(x,y) H\big(f(y)-f(x)\big),
\end{equation}
and 
\begin{equation}
\Upsilon(z)=e^z-1-z,\quad z\in \iR.
\end{equation}
Recall the classical case, where we have $\partial_t v-\Delta v=|\nabla v|^2=\Gamma(v)$.
In the nonlocal case, the carr\'e du champ operator is given by
\[
\Gamma(f)(x)=\frac{1}{2}\sum_{y\in X}k(x,y) \big(f(y)-f(x)\big)^2,
\]
which clearly differs from $\Psi_\Upsilon(f)$. Indeed, $\Gamma(f)$ corresponds to the first
term in the Taylor expansion of 
$\Upsilon(z)$ at zero as it appears in $\Psi_\Upsilon(f)$. In this sense, the function $\Upsilon$ in the nonlocal case plays the role of the square function in the local setting.
\section{The $CD_\Upsilon$-condition}
We again consider the Markov chain setting with finite or countably infinite state space $X$ and the generator given by  
\[
Lf(x)=\sum_{y\in X}k(x,y)\big(f(y)-f(x)\big),
\]
where $k(x,y)\ge 0$ for $x\neq y$, and $\sum_{y\in X}k(x,y)=0$ for all $x\in X$.
The question arises whether there is a natural analogue of $CD(\kappa,d)$ with corresponding implications, such as Li--Yau type inequalities and functional inequalities. Recall that $CD(\kappa,d)$ means that $\Gamma_2(f)\ge \kappa \Gamma(f)+\frac{1}{d}(Lf)^2$. Can one identify appropriate discrete/nonlocal counterparts of $\Gamma$, $\Gamma_2$ and the dimensional term $\frac{1}{d}(Lf)^2$?

In the diffusive setting, that is, $L$ satisfies the chain rule, the operators $\Gamma$ and $\Gamma_2$
can be well motivated by means of entropy, which is given by
\[
{\cH}(P_t f)=\int_X P_t f\log(P_t f)\,d\mu,
\]
where $P_t$ denotes the Markov semigroup generated by $L$ and $\mu$ is an invariant and reversible probability measure. In fact, if the chain rule is valid for $L$, we have
\begin{align*}
\frac{d}{dt}{\cH}(P_t f) &=-\int_X P_t f\, \Gamma\big(\log(P_t f)\big)\,d\mu, \\
\frac{d^2}{dt^2}{\cH}(P_t f) &=2\int_X P_t f \,\Gamma_2\big(\log(P_t f)\big)\,d\mu.
\end{align*}
While the first identity follows immediately, the second one requires additional manipulations, namely an integration-by-parts argument that relies on reversibility, see e.g.\ \cite{WZ}. Note that
$CD(\kappa,\infty)$ with $\kappa>0$ implies a differential inequality for the Fisher information,  defined by
\[
\cI(P_t f)=\int_X P_t f\, \Gamma\big(\log(P_t f)\big)\,d\mu.
\]
This observation constitutes the key idea of the Bakry--\'Emery approach to the logarithmic Sobolev inequality and the exponential decay of entropy.

It turns out that the identities for the first and second derivative of the entropy can be transferred to the discrete case. Weber and the author (\cite{WZ}) obtained the relations
\begin{align}
\frac{d}{dt}{\cH}(P_t f) &=-\int_X P_t f\, \Psi_\Upsilon\big(\log(P_t f)\big)\,d\mu,\label{EntDiss}\\
\frac{d^2}{dt^2}{\cH}(P_t f)&=2\int_X P_t f \,\Psi_{2,\Upsilon}\big(\log(P_t f)\big)\,d\mu,\nonumber
\end{align}
where $\Psi_\Upsilon$ is exactly the operator introduced in the previous section (see \eqref{PsiH}) and
\[
\Psi_{2,\Upsilon}(f)=\frac{1}{2}\big(L\Psi_\Upsilon(f)-B_{\Upsilon'}(f,Lf)\big),
\]
with ($\Upsilon'(z)=e^z-1$)
\[
B_{\Upsilon'}(f,g)(x)=\sum_{y\in X}k(x,y)\Upsilon'\big(f(y)-f(x)\big)\big(g(y)-g(x)\big).
\]
This motivated the following definition, see \cite{Web,WZ}.
\begin{definition} \label{Upsdef}
The operator $L$ is said to satisfy the condition $CD_\Upsilon(\kappa,F)$ with $\kappa\in \iR$ and CD-function $F$ if for all (suitable) $f:\,X\to \iR$ (and with $F_0:=F \chi_{[0,\infty)}$)
\[
\Psi_{2,\Upsilon}(f)\ge \kappa \Psi_\Upsilon(f)+F_0(-Lf)\quad \mbox{on}\;X.
\]
Further, $L$ is said to satisfy the condition $CD_\Upsilon(\kappa,d)$ with $\kappa\in \iR$ and $d\in [1,\infty]$ if
\[
\Psi_{2,\Upsilon}(f)\ge \kappa \Psi_\Upsilon(f)+\frac{1}{d}(Lf)^2\quad \mbox{on}\;X.
\]
\end{definition}
In analogy to the classical Bakry--\'Emery condition $CD(\kappa,d)$, the $CD_\Upsilon$ condition 
entails several important consequences, see \cite{WZ}. Firstly, $CD_\Upsilon(\kappa,\infty)$ with $\kappa>0$ implies the modified log-Sobolev inequality (MLSI) with constant $\kappa$, which means that
\[
{\cH}(f)\le \frac{1}{2\kappa}{\cI}(f)\quad \text{for all}\;f\in 
{\cP}(X).
\]
Here we set $d\mu=\pi d\#$, where $\#$ denotes the counting measure, and define
\[
{\cP}(X)=\big\{\rho:\,X\to [0,\infty)\;\mbox{s.t.}\;\int_X \rho \,d\mu=\sum_{x\in X}\rho(x)\pi(x)=1\big\}.
\]
The (relative) entropy is given by
\[
{\cH}(f)=\int_X f \log f\,d\mu=\sum_{x\in X} f(x) \log (f(x))\,\pi(x).
\]
Using detailed balance, i.e.\ $\pi(x)k(x,y)=\pi(y)k(y,x)$, the Fisher information can be written as
\begin{align*}
{\cI}(f) &= \int_X f \,\Psi_\Upsilon(\log f)\,d\mu\\
& = \frac{1}{2}\sum_{x,y\in X}k(x,y)\big(f(y)-f(x)\big) \big(\log f(y)-\log f(x)\big)\pi(x).
\end{align*}
Observe that \eqref{EntDiss} and the MLSI with constant $\kappa>0$ imply that
$\cH(P_t f)\le e^{-2\kappa t}\cH(f)$, $t\ge 0$, for all $f\in \cP(X)$.

Moreover, $CD_\Upsilon(\kappa,\infty)$ is characterised by the gradient bound
\[
\Psi_\Upsilon(P_t f)\le e^{-2\kappa t}P_t\big(\Psi_\Upsilon(f)\big),\quad t>0.
\]
Curvature bounds are also preserved under tensorisation. More precisely, let $L_i$ generate a Markov chain on $X_i$ and satisfy $CD_\Upsilon(\kappa_i,\infty)$ for $i=1,2$. Then the product generator $L:=L_1\oplus L_2$, acting on $X_1\times X_2$, satisfies $CD_\Upsilon(\kappa,\infty)$ with $\kappa=\min\{\kappa_1,\kappa_2\}$. Furthermore, the condition $CD_\Upsilon$ is compatible with the diffusion setting: there exists a natural notion of $\Gamma\oplus \Psi_\Upsilon$ and of $(\Gamma\oplus \Psi_\Upsilon)_2$, for which an analogous tensorisation principle holds for hybrid processes. Finally, we note that $CD_\Upsilon(\kappa,d)$ implies $CD(\kappa,d)$; see \cite{WZ} for these results as well as for a detailed discussion of the relation between the $CD_\Upsilon$ condition and other curvature-dimension conditions in the literature.

To illustrate these concepts, consider first the unweighted complete graph $K_n$, i.e.\ the graph in which every pair of the $n$ vertices is connected. In this case, $CD_\Upsilon(\sqrt{2n},\infty)$ holds for all $n\ge 2$. This bound is optimal for $n=2$, since $CD(2,\infty)$ is sharp in that case. More generally, $CD(1+\frac{n}{2},\infty)$ is optimal for all $n\ge 2$. As a second example, the $n$-dimensional hypercube satisfies $CD_\Upsilon(2,\infty)$ for all $n\in \mathbb{N}$, and this bound is optimal, as $CD(2,\infty)$ is sharp for every $n$. This can be seen by induction, using the result for $K_2$ and the tensorisation property. Finally, the unweighted $3$-star fails $CD_\Upsilon(0,\infty)$ at the center point, but $CD_\Upsilon(\kappa,\infty)$ with some $\kappa\in (-\infty,0)$ is true.
\section{Long-range jump operators}
We consider now nonlocal operators on the lattice $\iZ$ of the form
\begin{equation} \label{genlaplacedef}
L f(x)=\,\sum_{j\in \iZ} k(j) \big(f(x+j)-f(x)\big),\quad x\in \iZ,
\end{equation}
with a kernel $k:\,\iZ\to [0,\infty)$ satisfying $0<\sum_{j\in \iZ}k(j)<\infty$, $k(-j)=k(j)$ for all $j\in \iN$, and $k(0)=0$. An important example is given by
\[
k_\beta(j)=\,\frac{c}{|j|^{1+\beta}},\quad j\in \iZ\setminus \{0\},\quad \mbox{with}\;c,\,\beta>0.
\]
For $\beta\in (0,2)$, the corresponding operator $L$ is closely linked to the fractional discrete Laplacian $-(-\Delta)^{\beta/2}$ on $\iZ$. The following statements, due to Spener, Weber and the author of these notes \cite{SWZ1}, address the validity of the classical CD-condition.
\begin{theorem} Let the kernel $k$ be as described before.
\begin{itemize}
\item [(i)] If $k$ is nonincreasing on $\iN$ and such that $\sum_{j\in \iN} k(j)j^2<\infty $, then $CD(0,d)$ holds for some $d<\infty$.
\item [(ii)] In particular, for $k=k_\beta$ with $\beta>2$,  we have $CD(0,d)$ with some 
$d<\infty$.
\item [(iii)] If $k=k_\beta$ and $\beta<2$, then $CD(0,d)$ fails for all $d<\infty$.
\end{itemize}
\end{theorem}
The third statement shows the limitations of the classical CD-condition. $CD_\Upsilon(0,F)$ turns out to be more flexible. What about $\Psi_{2,\Upsilon}(f)\ge F_0(-Lf)$ with some $F_0$ as in Definition \ref{Upsdef}? Note that 
\[
\Psi_{2,\Upsilon}(f) (x) = \frac{1}{2} \sum\limits_{j,l \in \mathbb{Z}} k(j) k(l) e^{f(x+l)-f(x)} \Upsilon\big( f(x+j+l)-f(x+j)-f(x+l)+f(x)\big),
\]
where $\Upsilon(z)=e^z-1-z$. The following two theorems were established by Kr\"{a}ss, Weber and the author of these notes \cite{KWZ}.
\begin{theorem} \label{CDUps2} Let the kernel $k$ be as described before.
\begin{itemize}
\item [(i)]  Suppose that $\sum_{j\in \iN} k(j)^{1-\delta}<\infty$ with some $\delta\in (0,1)$. 
Then 
\[
\Psi_{2,\Upsilon}(f)\ge c|Lf|^\gamma,\quad\quad  \gamma=\frac{1+\delta}{\delta}.
\]
\item [(ii)]  If $k=k_\beta$ with $\beta\in (0,\infty)$, then $CD_\Upsilon(0,F)$ holds with some CD-function $F$, which grows exponentially at $\infty$ and satisfies $F(x)\sim cx^\gamma$ as $x\to 0$, where $\gamma=2$ for $\beta>2$ and  ($\beta_*:=
\frac{1+\sqrt{5}}{2}$)
\[
\gamma=\frac{1+2\beta}{\beta}\quad\mbox{for}\;\,\beta\in (0,\beta_*], \quad\gamma=\frac{\beta-\varepsilon}{\beta-\varepsilon-1}\quad\mbox{for}\;\,\beta\in (\beta_*,2],
\]
with arbitrary small $\varepsilon>0$.
\end{itemize}
\end{theorem}
With these positive results at hand, we can derive Li--Yau and Harnack inequalities for the heat equation with jump operator $L$.
\begin{theorem} \label{CDUps3}
Let $\beta \in (0,\infty)$ and $L_\beta$ be the jump operator associated with $k_\beta$. Then every bounded function $u : [0,\infty) \times \iZ \to (0,\infty)$ that is $C^1$ in time and solves $\partial_t u = L_\beta u$
 on $ (0,\infty)\times \iZ$ satisfies the Li--Yau estimate
\[
-L_\beta \log u (t,x) \leq \varphi(t),\quad (t,x) \in (0,\infty)\times \iZ,
\]
where $\varphi$ is the relaxation function corresponding to $F$ from Theorem \ref{CDUps2} above. 
We have $\varphi(t) \sim -c\log t$ as $t\to 0$, and $\varphi(t) \sim c t^{-\frac{1}{\gamma-1}}$ as $t\to \infty$.

Moreover, we have the 
Harnack inequality
\begin{equation*}
u(t_1,x_1) \leq u(t_2,x_2) \exp \Big(\int_{t_1}^{t_2} \varphi(t) \ \mathrm{d}t + \frac{2 \vert x_1-x_2 \vert^{\min\{1+\beta,2\}}}{t_2-t_1}\Big),
\end{equation*}
for $0 \leq t_1 < t_2$, $x_1,x_2 \in \mathbb{Z}$.
\end{theorem}
Note that, in contrast to the classical parabolic case, $t_1=0$ is admissible. This is due to the integrability of
$\varphi$ at $0$. Moreover, Theorem \ref{CDUps2} (ii) yields, for $\beta\in (0,2)$, an analogous result for the fractional {\em discrete} Laplacian $-(-\Delta)^{\frac{\beta}{2}}$ on $\iZ$. 

Can we derive similar results for the fractional Laplacian in the continuous setting? 
Recall that for $\beta\in (0,2)$, the operator  $-(-\Delta )^\frac {\beta}{2}$ on $\iR^d$ is given by \eqref{frakLaplacianDef}, where $c_{\beta,d}$ is a normalising constant. 
The corresponding carr\'{e} du champ and iterated carr\'{e} du champ take the form
\begin{align*}                                                                                                         
\Gamma(u)(x) & = {c_{\beta,d}}\int_{\iR^d}\frac{(u(x+h)-u(x))^2}{|h|^{d+\beta}}\, dh, \\
\Gamma_2(u)(x) & =  {c_{\beta,d}^2}\!\int_{\iR^d}\!\!\int_{\iR^d}\!\!\!\frac{[u(x\!+\!h+\!\sigma)\!-\!u(x\!+\!h)\!-\!u(x\!+\!\sigma)\!+\!u(x)]^2}{|h|^{d+\beta}|\sigma|^{d+\beta}} dh\, d\sigma.
 \end{align*}
The following result, due to Spener, Weber, Zacher \cite{SWZ2}, gives a negative answer with regard to the approach via CD inequalities.
\begin{theorem}
For any $R > 0$, $\kappa \in \iR$ and $N\in (0,\infty)$ there is $u \in {C}^\infty_c(\iR^d)$ such that
\begin{equation*}
 0 < \Gamma_2(u)(x) < \kappa \Gamma(u)(x) + \frac{1}{N} \left({L}(u)(x)\right)^2,\quad \forall x\in B(0,R).
\end{equation*}
\end{theorem}
Note that, since $CD(0,N)$ fails, $CD_\Upsilon(0,N)$ must also fail for all finite $N\ge 1$.
Thus, the approach based on CD inequalities offers little hope. Nevertheless, a Li--Yau inequality for the fractional Laplace operator can be proved using a different method based on the associated heat kernel.
\section{Reduction to the heat kernel}
Let us first explain the reduction method in the classical case of the Laplacian. A key ingredient is the following observation. Let $Pf(x)=\int_{\iR^d}H(x,y)f(y)\,dy$, for sufficiently regular, positive functions $H$ and $f$. Then
\begin{equation} \label{Unglclassic}
\int_{\iR^d}\big|\nabla_x \log H(x,y)\big|^2 H(x,y) f(y)\,{d}y\ge \big|\nabla \log Pf(x)\big|^2 Pf(x),\quad x\in \iR^d.
\end{equation}
Indeed, by H\"older's inequality we have, for $i=1,\ldots,d$,
\begin{align*}
(\partial_{x_i} Pf(x))^2 & = \big(\int_{\iR^d} \partial_{x_i} H(x,y)f(y)\,{d}y\big)^2\\
& \le \int_{\iR^d} \frac{(\partial_{x_i} H(x,y))^2}{H(x,y)}f(y) \,{d}y \,\,
\int_{\iR^d}  H(x,y) f(y) \,{d}y,
\end{align*}
which directly leads to \eqref{Unglclassic} by summing up and employing the chain rule for the gradient 
($\nabla (\log g)=\frac{\nabla g}{g}$).

For the heat kernel  $H(t,x,y)=(4\pi t)^{-\frac{d}{2}}e^{-\frac{|x-y|^2}{4t}}$, we have
\begin{equation} \label{HK1}
-\Delta_x \big(\log H(t,x,y)\big)=\frac{d}{2t}=:\varphi(t).
\end{equation}
For any positive (smooth) solution $u$ of the heat equation,
\[
\partial_t u-u\Delta(\log u)=u|\nabla (\log u)|^2.
\]
In particular,
\[
\partial_t H(t,x,y)+H(t,x,y)\varphi(t) =H(t,x,y)|\nabla_x (\log H(t,x,y))|^2.
\]
Consider a positive solution $u(t,x)=\int_{\iR^d} H(t,x,y) u_0(y)\,dy$ of the heat 
equation (cf.\ the Widder theorem). Then
\begin{align*}
\partial_t u(t,x)&+\varphi(t)u(t,x)=\int_{\iR^d}\big(\partial_t H(t,x,y) +\varphi(t)H(t,x,y)\big)u_0(y)\,dy\\
& = \int_{\iR^d}\big(  |\nabla_x (\log H(t,x,y))|^2  H(t,x,y)\big)u_0(y)\,dy\\
& \ge |\nabla(\log u(t,x))|^2 u(t,x)\quad\quad\quad \text{(by \eqref{Unglclassic})}\\
& = \partial_t u(t,x)-u(t,x)\Delta(\log u(t,x)).
\end{align*}
Hence
\[
-\Delta(\log u(t,x))\le \varphi(t)=\frac{d}{2t}.
\]
Note that here we only need $-\Delta_x \big(\log H(t,x,y)\big)\le \varphi(t)$. This argument shows that the Li--Yau inequality for the heat kernel implies the same Li--Yau estimate for positive solutions.

Is there a similar argument for the fractional heat equation (FHE)?
Let $u$ be a positive solution of the FHE. Then, analogously to \eqref{loguEqu} and with $\beta\in(0,2)$,
\[
\partial_t u+u(-\Delta)^{\frac{\beta}{2}}(\log u)=u\Psi_\Upsilon(\log u),
\]
where
\[
\Psi_\Upsilon(v)(x)=c_{\beta,d}\int_{\R^d} \frac{\Upsilon\big(v(y)-v(x)\big)}{|x-y|^{d+\beta}}\,dy.
\]
It turns out that there exists a nonlocal analogue of the key inequality 
\eqref{Unglclassic}. Letting $P, H, f$ be as before, one has
\begin{equation} \label{keyUngl}
\int_{\iR^d} \Psi_\Upsilon (\log H(\cdot,y))(x)H(x,y)f(y)\,dy \geq \Psi_\Upsilon(\log P f)(x) P f(x),\quad x\in \iR^d,
\end{equation}
see \cite[Lemma 2.2]{WZMA}. Moreover, it is known from \cite{BPSV} that positive strong solutions $u$ of the FHE (see \cite{BPSV} for the precise definition) admit the representation
\begin{equation}\label{eq:solutionfractionalheatequation}
u(t,x) = \int_{\R^d}G^{(\beta)}(t,x-y) u_0(y)\,{d}y,
\end{equation}
where $G^{(\beta)}$ is the fundamental solution of the FHE. Using \eqref{keyUngl}, one can argue as before to obtain the implication
\[
(-\Delta)^{\frac{\beta}{2}}(\log G^{(\beta)} )(t,x)\le \varphi(t)\quad \Rightarrow \quad 
(-\Delta)^{\frac{\beta}{2}}(\log u)(t,x)\le \varphi(t).
\]
This raises the question for which $\varphi$ the left-hand side inequality holds. In \cite[Lemma 3.1]{WZMA}, Weber and Zacher show that there exists a constant $C_{LY}>0$, depending only on $\beta$ and $d$, such that
for all $\beta \in (0,2)$,
\[
(-\Delta )^{\frac{\beta}{2}}(\log G^{(\beta)})(t,x) \leq \frac{C_{LY}}{t}, \quad t>0,\,x\in \iR^d.
\]
With this inequality at hand, they obtain the following result, see \cite[Theorem 3.2]{WZMA}.
\begin{theorem}
Let $\beta\in (0,2)$, and let $u:\,[0,\infty)\times \R^d \to (0,\infty)$ be a strong solution of
\[
\partial_t u+(-\Delta)^{\frac{\beta}{2}}u=0\quad \mbox{in}\;(0,\infty)\times \iR^d.
\]
Then for all $(t,x)\in (0,\infty)\times \R^d$, $u$ satisfies the Li--Yau inequality
\begin{equation}\label{eq:fractionalLiYau}
\big(-\Delta\big)^\frac{\beta}{2} (\log u)(t,x) \leq \frac{C_{LY}}{t}.
\end{equation}
Moreover, this is equivalent to the differential Harnack inequality
\begin{equation}\label{eq:fractionalDiffHar}
\partial_t(\log u)(t,x) \geq \Psi_\Upsilon(\log u)(t,x)-\frac{C_{LY}}{t}.
\end{equation}
\end{theorem}

Weber and Zacher also show that the differential Harnack inequality \eqref{eq:fractionalDiffHar}
implies the Harnack inequality
\begin{equation} \label{HarnackFL}
u(t_1,x_1)\le u(t_2,x_2) \Big(\frac{t_2}{t_1}\Big)^{C_{LY}}\exp \left(C\left[1+\frac{|x_1-x_2|^{\beta+d}}
{(t_2-t_1)^{1+\frac{d}{\beta}}}\right]\right),
\end{equation}
for all $0<t_1<t_2<\infty$ and $x_1,x_2 \in \mathbb{R}^d$, where the constant $C>0$ only depends on $\beta$ and $d$, see \cite[Theorem 5.2]{WZMA}. Note that parabolic Harnack inequalities for the space fractional heat equation were already known before, even in more general settings, see e.g.
\cite{BL2,BoSiVa,CLD,CK} and the more recent work \cite{KaWe}. The key point here is that the differential Harnack inequality \eqref{eq:fractionalDiffHar} is sufficiently strong to 
yield the parabolic Harnack inequality, thereby also justifying its name.

In the case of the FHE, in contrast to the classical heat equation, there also exist backward and elliptic versions of the parabolic Harnack inequality, that is, the condition $0<t_1<t_2<\infty$ is not required, cf.\ \cite{BoSiVa}. The reason is that the fractional heat kernel $G^{(\beta)}$, for $\beta\in (0,2)$, satisfies 
\begin{equation} \label{ExtraEst}
t|\partial_t G^{(\beta)}(t,x)|\le c_{\beta,d}G^{(\beta)}(t,x),\quad t>0,\,x\in \iR^d.
\end{equation}
Using \eqref{ExtraEst}, the argument of Weber and Zacher in \cite{WZMA} can be refined to yield a differential Harnack estimate of the form
\begin{equation} \label{betterDH}
|\partial_t(\log u)| + \Psi_\Upsilon(\log u)\le \frac{\tilde{C}(\beta,d)}{t},\quad t>0,\,x\in \iR^d.
\end{equation}
On the basis of \eqref{betterDH}, the proof of \cite[Theorem 5.2]{WZMA} can then be adapted to also cover Harnack estimates of backward and elliptic type. 


\begin{thebibliography}{99}
{\footnotesize
\bibitem{BGL} 
Bakry, D.; Gentil, I.; Ledoux, M.: \textit{Analysis and geometry of Markov diffusion operators}. Grundlehren der mathematischen Wissenschaften \textbf{348}. Springer Cham, Heidelberg, 2014.
\bibitem{BE}
Bakry, D.; \'Emery, M.: \textit{Diffusions hypercontractives}. Séminaire de probabilités, XIX, 1983/84, 177-206. Lecture Notes in Math., Springer, Berlin, 1985.
\bibitem{BPSV}
Barrios, B.; Peral, I.; Soria, F.; Valdinoci, E.:
A Widder's type theorem for the heat equation with nonlocal diffusion. 
Arch. Ration. Mech. Anal. {\bf 213} (2014), no. 2, 629--650.
\bibitem{BL2}
Bass, R.\ F.; Levin, D.\ A.: Transition probabilities for symmetric jump processes. Trans.\ Amer.\ Math.\ Soc.\ {\bf 354}
(2002), 2933--2953.
\bibitem{BoSiVa}
Bonforte, M.; Sire, Y.; V\'azquez, J.\ L.: Optimal existence and uniqueness theory for the fractional heat equation. Nonlinear Anal. {\bf 153} (2017), 142--168.
\bibitem{CLD}
Chang-Lara, H.\ A.; D'avila, G.: H\"older estimates for non-local parabolic equations with critical drift. J. Differential
Equations {\bf 260} (2016), 4237--4284.
\bibitem{CK}
Chen, Z.-Q.; Kumagai, T.: Heat kernel estimates for stable-like processes on d-sets. Stochastic Process. Appl.
{\bf 108} (2003), 27--62.
\bibitem{DKZ} 
Dier, D.; Kassmann, M.; Zacher, R.: Discrete versions of the Li--Yau gradient estimate. 
Ann.\ Sc.\ Norm.\ Super.\ Pisa Cl.\ Sci.\ (5) {\bf 22} (2021), 691--744.
\bibitem{EM}
Erbar, M.; Maas, J.:
Ricci curvature of finite Markov chains via convexity of the entropy.
Arch.\ Ration.\ Mech.\ Anal.\ {\bf 206} (2012), 997--1038. 
\bibitem{Harvard} 
Bauer, F.; Horn, P.; Lin, Y.; Lippner, G.; Mangoubi, D.; Yau, S.-T.: Li--Yau inequality on graphs. J. Differential Geom. \textbf{99} (2015), 359--405.
\bibitem{Garo}
Garofalo, N.:
Fractional thoughts.
Danielli, D. (ed.) et al., New developments in the analysis of nonlocal operators. AMS special session, University of St. Thomas, Minneapolis, MN, USA, October 28–30, 2016. Providence, RI: American Mathematical Society (AMS). Contemp.\ Math.\ {\bf 723} (2019), 1--135. 
\bibitem{KaWe}
Kassmann, M.; Weidner, M.:
The parabolic Harnack inequality for nonlocal equations.
Duke Math. J. {\bf 173} (2024), 3413--3451. 
\bibitem{KWZ}
Kr\"{a}ss, S.; Weber, F.; Zacher, R.:
Li--Yau and Harnack inequalities via curvature-dimension conditions for discrete long-range jump operators including the fractional discrete Laplacian.
Discrete Contin. Dyn. Syst. {\bf 44} (2024), 1982--2028. 
\bibitem{LY} 
Li, P.; Yau, S.-T.: On the parabolic kernel of the Schr\"odinger operator. Acta Math. {\bf 156} (1986), 153--201.
\bibitem{Mie}
Mielke, A.:
Geodesic convexity of the relative entropy in reversible Markov chains.
Calc.\ Var.\ Partial Differ.\ Equ.\ {\bf 48} (2013), 1--31. 
\bibitem{MUN} 
M\"unch, F.: Li--Yau inequality on finite graphs via non-linear curvature dimension conditions. J.\ Math.\ Pures Appl. (9) {\bf 120} (2018), 130--164. 
\bibitem{Olli}
Ollivier, Y.:
Ricci curvature of Markov chains on metric spaces.
J.\ Funct.\ Anal.\ {\bf 256} (2009), 810--864. 
\bibitem{SWZ1}
Spener, A.; Weber, F.; Zacher, R.:
Curvature-dimension inequalities for non-local operators in the discrete setting. Calc.\ Var.\ Partial Differ.\ Equ. {\bf 58} (2019), Paper No. 171. 
\bibitem{SWZ2}
Spener, A.; Weber, F.; Zacher, R.:
The fractional Laplacian has infinite dimension. Commun.\ Partial Differ.\ Equations {\bf 45} (2020), 57--75. 
\bibitem{Web}
Weber, F.:
Entropy-information inequalities under curvature-dimension conditions for continuous-time Markov chains.
Electron.\ J.\ Probab.\ {\bf 26} (2021), Paper No. 52, 31 p. 
\bibitem{WZMA}
Weber, F.; Zacher, R.:
Li--Yau inequalities for general non-local diffusion equations via reduction to the heat kernel.
Math.\ Ann.\ {\bf 385} (2023), 393--419. 
\bibitem{WZ}
Weber, F.; Zacher, R.:
The entropy method under curvature-dimension conditions in the spirit of Bakry--\'Emery in the discrete setting of Markov chains.
J.\ Funct.\ Anal.\ {\bf 281} (2021), Article ID 109061, 81 p. 
}
\end{thebibliography}
\end{document}